\documentclass[12pt]{aptpub}
\renewcommand{\crnotice}{\hfill{{\small\em
\/}\ {\small}}}

\usepackage{amsmath,amstext,url}
\usepackage{color}
\usepackage{amsfonts}
\usepackage{amsmath}
\usepackage{txfonts}
\usepackage{mathrsfs}
\usepackage{algorithm}
\usepackage{placeins}

\renewcommand{\authornames}{\underline{\footnotesize }
J. Li $\&$ M. Hou}
\renewcommand{\shorttitle}
{\underline{Long-time behaviour of GWSs with circular mechanism}}

\numberwithin{equation}{section} 

\makeatletter \@addtoreset{equation}{section}

\makeatletter \@addtoreset{lemma}{section}

\makeatletter \@addtoreset{theorem}{section}

\makeatletter \@addtoreset{proposition}{section}

\makeatletter \@addtoreset{corollary}{section}

\makeatletter \@addtoreset{remark}{section}

\makeatletter \@addtoreset{definition}{section}

\makeatletter \@addtoreset{example}{section}

\makeatletter \@addtoreset{algorithm}{section}

\begin{document}

\thispagestyle{firstpg}

\vspace*{1.2pc} \noindent \normalsize\textbf{\Large {Long-time behaviour of Galton-Watson systems with circular mechanism}} \hfill

\vspace{12pt} \hspace*{0.75pc}{\small\textrm{Junping Li}}
\hspace{-2pt}$^{*}$, {\small\textit{Guangdong University of Science and Technology; Central
South University}}

\vspace{5pt}
\hspace*{0.75pc}{\small\textrm{Mixuan Hou}}
\hspace{-2pt}$^{**}$, {\small\textit{Central South University}}
\par
\footnote{\hspace*{-0.75pc}$^{*}\,$Corresponding author. Postal
address: Guangdong University of Science and Technology, Dongguan, 523083, China; Central
South University, Changsha, 410083, China. E-mail:
jpli@mail.csu.edu.cn}
\footnote{\hspace*{-0.75pc}$^{**}\,$Postal address:
 School of Mathematics and Statistics, Central South University, Changsha, 410083, China. E-mail address:
242101015@csu.edu.cn}

\par
\par
\renewenvironment{abstract}{%
\vspace{8pt} \vspace{0.1pc} \hspace*{0.25pc}
\begin{minipage}{14cm}
\footnotesize
{\bf Abstract}\\[1ex]
\hspace*{0.5pc}} {\end{minipage}}
\begin{abstract}
Let $\{Z_n:n\geq 0\}$ be a Galton-Watson system with circular mechanism $\emph{\textbf{a}}*\emph{\textbf{b}}$, where $\emph{\textbf{a}}=\{a_j\}_{j=0}^{\infty}$ and $\emph{\textbf{b}}=\{b_j\}_{j=0}^{\infty}$ are probability distributions on $\mathbf{Z}_+:=\{0,1,2,\cdots\}$. Let $m_a:=\sum\limits_{j=0}^{\infty}ja_j$,  $m_b:=\sum\limits_{j=0}^{\infty}jb_j$. The extinction property of such branching systems is first studied. Then, it is proved that there exists $\Gamma_n$ such that $W_n=\Gamma_n^{-1} Z_n$ is an integrable martingale and hence converges to some random variable $W$. Moreover, for the case that $a_0=b_0=0, a_1,b_1\in (0,1)$, the convergence rates to $0$ of
\begin{eqnarray*}
P(|\frac{Z_{n+1}}
{Z_n}-\omega_n|>\varepsilon\mid Z_0=1),\quad
 P(|W_n-W| >\varepsilon\mid Z_0=1)
\end{eqnarray*}
and
\begin{eqnarray*}
P(|\frac{Z_{n+1}}
{Z_n}-\omega_n|>\varepsilon\mid W\geq \delta, Z_0=1)
\end{eqnarray*}
as $n\to \infty$ are presented for $\varepsilon,\delta>0$ under various moment conditions on $\{a_j\}$ and $\{b_j\}$, where $\omega_n=m_a$ or $m_b$ for $n$ being even or odd respectively. It is further shown that the first rate is geometric while the last two rates are supergeometric under a finite moment generating function hypothesis.
\end{abstract}

\vspace*{12pt}
\parbox[b]{26.75pc}{{
}}
\par

\vspace*{12pt}
{\footnotesize {\bf \ MSC:} 60J80


\par
\vspace{5mm}
 \setcounter{section}{1}
 \setcounter{equation}{0}
 \setcounter{theorem}{0}
 \setcounter{lemma}{0}
 \setcounter{corollary}{0}

\noindent {\large \bf 1. Introduction}
\vspace{3mm}
\par
The Markov branching systems (MBSs) play an important role in the classical field of stochastic systems. The basic property of MBS is the branching property, i.e., different individuals act independently when giving birth or death and the system stops when there is no particle in it. The classical Markov branching systems are very deeply studied, standard references are Asmussen $\&$ Jagers~\cite{AJ97}, Asmussen $\&$ Hering~\cite{ASHH83}, Athreya $\&$ Ney~\cite{AKN72} and Harris~\cite{Harr63}. Based on the standard branching structure, some generalized branching systems are well studied. For example, Yamazato~\cite{YM75} investigated a branching systems with immigration which only occurs at state zero. Chen~\cite{CA2002} considered general branching systems with or without resurrection. Sevast'yanov~\cite{SBA49} and Vatutin~\cite{VVA74} considered the interacting branching systems. Liu $\&$ Zhang~\cite{LZ16} studied large deviation for supercritical branching processes with immigration. Zhang, Li $\&$ Geng~\cite{ZLG2024} discussed nonlinear Markov branching processes with immigration and resurrection. Mitov $\&$ Yanev~\cite{MY2025} considered a class of critical Markov branching processes with non-homogeneous Poisson immigration. Hermann $\&$ Pfaffelhuber~\cite{HP2020} investigated the extinction, survival and duality for Markov branching processes with disasters. Furthermore, Pasha~\cite{P2020} considered the stability of traveling wave solutions for integro-differential equations related to branching Markov processes. Pakes~\cite{PA2025} studied the stationary measure for Markov branching processes. Imomov $\&$ Murtazaev~\cite{IM2024} discussed the Kolmogorov constant explicit form for discrete-time stochastic branching systems. Francisci $\&$ Vidyashnkar~\cite{FV2024} studied branching processes in random environments with Thresholds. Mitov $\&$ Yanev~\cite{MY2024} considered critical Markov branching process with infinite variance allowing Poisson immigration with increasing intensity.  Simon, Emma, Andreas $\&$ Ellen~\cite{SEAE2024} investigated the properties of non-local branching Markov process. Smorodina $\&$ Yarovaya~\cite{SY2024} studied the limit behaviour of branching random walks.
\par
It is well-known that the evolution of a branching system is controlled by its branching mechanism. However, in realistic situations, such as controlled population models and controlled molecular biology models, the controller may change the branching mechanism at different time or states. Therefore, the evolution behavior of the system will be controlled by all the branching mechanisms involved.
\par
In order to clearly describe the model considered in this paper, we first give the following definitions.
\par
Let $(\Omega, \mathscr{F}, P)$ be a probability space and denote
\begin{eqnarray}\label{eq1-1}
\mathscr{P}=\{\emph{\textbf{a}}=\{a_k\}_{k=0}^{\infty}:\ a_k\geq 0, \forall k\geq 0\ {\rm{and}}\ \sum_{k=0}^{\infty}a_k=1\},
\end{eqnarray}
i.e., $\mathscr{P}$ is the set of all probability distributions on $\mathbf{Z}_+:=\{0,1,2,\cdots\}$. Obviously, $\mathscr{P}$ is a Borel subset of Banach space $l_{\infty}$. An element $\emph{\textbf{a}}\in \mathscr{P}$ is also called a branching mechanism in a branching system, which is the offspring distribution of the particles in the system. For any $\emph{\textbf{a}}=\{a_k\}_{k=0}^{\infty}\in \mathscr{P}$, define
\begin{eqnarray*}
  f(\emph{\textbf{a}};s)=\sum_{k=0}^{\infty}a_ks^k,\ \ |s|\leq 1.
\end{eqnarray*}
Denote $m_a=f'(\emph{\textbf{a}};1)=\sum_{k=0}^{\infty}ka_k$ and $\rho_a$ the smallest nonnegative root of $f(\emph{\textbf{a}};s)=s$.
\par
\begin{definition}\label{def1.1}
Let $\{\emph{\textbf{a}}^{(n)}\}_{n=0}^{\infty}$ be a sequence of branching mechanisms in $\mathscr{P}$.
\par
(i)\ A $\mathbf{Z}_+$-valued system $\{Z_n: n\geq 0\}$ is called a Galton-Watson system in deterministic environment $\{\emph{\textbf{a}}^{(n)}\}_{n=0}^{\infty}$, if
\begin{eqnarray*}
Z_{n+1}=\sum\limits_{i=1}^{Z_n}\xi_{n,i}^{(\emph{\textbf{a}}^{(n)})},\ n\geq 0
\end{eqnarray*}
where $\{\xi_{n,i}^{(\emph{\textbf{a}}^{(n)})}:n\geq 0, i\geq 1\}$ is a sequence of independent $\mathbf{Z}_+$-valued random variables satisfying $E[s^{\xi_{n,i}^{(\emph{\textbf{a}}^{(n)})}}
]=f(\emph{\textbf{a}}^{(n)};s)$ for $i\geq 1$.
\par
(ii)\ If there exists an integral $d$ such that $\emph{\textbf{a}}^{(n)}=\emph{\textbf{a}}^{(k)}\ {\rm{if}}\ n=k\ mod(d+1)$, then $\{Z_n:n\geq 0\}$ is called a branching system with circular mechanism, or is simply called an $\emph{\textbf{a}}^{(0)}\ast \cdots \ast \emph{\textbf{a}}^{(d)}$-Galton-Watson system.
\end{definition}
\par
If $\emph{\textbf{a}}^{(n)}=\emph{\textbf{a}}\in \mathscr{P}$, then the process defined in Definition~\ref{def1.1} becomes the standard Galton -Watson process. In such case,  Athreya~\cite{AK94} considered the above convergence rates for supercritical Galton-Watson systems. Based on Athreya \cite{AK94}. Liu $\&$ Zhang~\cite{LZ16} studied the convergence rate of (\ref{eq1-4}) for Galton-Watson systems with immigration. Li $\&$ Li~\cite{L1} discussed the convergence rates of (\ref{eq1-5}) and (\ref{eq1-6}) for Galton-Watson systems with immigration and proved that they are supergeometric. Li, Cheng $\&$ Li~\cite{L2} investigated the above convergence rates for single-type continuous time branching systems.
\par
 Recall that the branching mechanism involved in the above references is time-independent. Motivated by the fact that the branching mechanism may be different at different times in realistic situations, in this paper, we mainly consider the long-time behaviour of $\emph{\textbf{a}}^{(0)}\ast \cdots \ast \emph{\textbf{a}}^{(d)}$-Galton-Watson systems. For simplicity, we only consider the case of $\emph{\textbf{a}}\!*\!\emph{\textbf{b}}$-Galton-Watson systems, where $\emph{\textbf{a}}, \emph{\textbf{b}}\in \mathscr{P}$. The general case can be similarly discussed.
\par
Let $\{Z_n:n\geq 0\}$ be an $\emph{\textbf{a}}\!*\!\emph{\textbf{b}}$-Galton-Watson system.
Denote $m_a:=f'(\emph{\textbf{a}};1),\ m_b:=f'(\emph{\textbf{b}};1)$ and
\begin{eqnarray}\label{eq1-2}
m:=m_am_b,\quad  \Gamma_n:=\begin{cases}
m^{\frac{n}{2}}, & \mbox{if}\ n\ \mbox{is\ even},\\
m^{\frac{n-1}{2}}m_a, & \mbox{if}\ n\ \mbox{is\ odd},
\end{cases}, \quad \omega_n=\begin{cases}
m_a, & \mbox{if}\ n\ \mbox{is\ even},\\
m_b, & \mbox{if}\ n\ \mbox{is\ odd}.
\end{cases}
\end{eqnarray}
Then, define
\begin{eqnarray}\label{eq1-3}
 W_n=\Gamma_n^{-1}Z_n,\ \ n\geq 0.
\end{eqnarray}
\par
More specifically, the main aim of this paper is to discuss the extinction property of $\emph{\textbf{a}}\!*\!\emph{\textbf{b}}$-Galton-Watson systems and the convergence rates of
\begin{eqnarray}\label{eq1-4}
P(|\frac{Z_{n+1}}
{Z_n}-\omega_n|>\varepsilon\mid Z_0=1),
\end{eqnarray}
\begin{eqnarray}\label{eq1-5}
 P(|W_n-W| >\varepsilon\mid Z_0=1),
\end{eqnarray}
and
\begin{eqnarray}\label{eq1-6}
P(|\frac{Z_{n+1}}
{Z_n}-\omega_n|>\varepsilon\mid W\geq \delta, Z_0=1)
\end{eqnarray}
as $n\to\infty$ for $\varepsilon, \delta>0$.

\vspace{5mm}
 \setcounter{section}{2}
 \setcounter{equation}{0}
 \setcounter{theorem}{0}
 \setcounter{lemma}{0}
 \setcounter{corollary}{0}
 \setcounter{assumption}{0}

\noindent {\large \bf 2. Extinction property of $\emph{\textbf{a}}\!*\!\emph{\textbf{b}}$-Galton-Watson systems}
\vspace{3mm}
\par
In this section, we discuss the basic property and the extinction behaviour of $\emph{\textbf{a}}\!*\!\emph{\textbf{b}}$-Galton-Watson systems. We first give some preliminaries.
\par
The following lemma~\ref{le2.1} is due to Athreya and Ney~\cite{AKN72} and the proof is omitted.
\par
\begin{lemma}\label{le2.1}For any $\textbf{a}=\{a_k\}_{k=0}^{\infty}\in \mathscr{P}$,
\ $f(\textbf{a};s)$ is a convex increasing function on $[0,1]$. If $m_a\leq 1$ then $f(\textbf{a};s)>s$ for
all $s \in [0,1)$ and $f(\textbf{a};s)=s$ has exactly one root $1$ on $[0,1]$. Furthermore, if
 $m_a<1$ then $1$ is a simple root while if $m_a=1$, then $1$ is a root of multiplicity $2$. If $m_a>1$, then $f(\textbf{a};s)=s$ has exactly two roots $\rho_a$ and $1$ on $[0,1]$ with $0\leq\rho_a<1$ such that $f(\textbf{a};s)>s$ for $s\in [0,\rho_a)$ and $f(\textbf{a};s)<s$ for $s\in (\rho_a,1)$. Both $\rho_a$
and $1$ are simple.
\end{lemma}
\par
For
$\emph{\textbf{a}}=\{a_k\}_{k=0}^{\infty}, \emph{\textbf{b}}=\{b_k\}_{k=0}^{\infty} \in \mathscr{P}$,
define
\begin{eqnarray*}
\alpha_0(s)=s,\ \alpha_1(s)=\alpha(s):=f(\emph{\textbf{a}};f(\emph{\textbf{b}};s)),\ \alpha_{n+1}(s)=\alpha(\alpha_n(s)),\ n\geq 0,
\end{eqnarray*}
and
\begin{eqnarray*}
\beta_0(s)=s,\ \beta_1(s)=\beta(s):=f(\emph{\textbf{b}};f(\emph{\textbf{a}};s)),\ \beta_{n+1}(s)=\beta(\beta_n(s)),\ n\geq 0.
\end{eqnarray*}
Let $\rho_a$, $\rho_b$, $\rho_{ab}$ and $\rho_{ba}$ denote the smallest nonnegative roots of $f(\emph{\textbf{a}};s)=s$, $f(\emph{\textbf{b}};s)=s$, $\alpha(s)=s$ and $\beta(s)=s$, respectively.
\par
\begin{lemma}\label{le2.2}
{\rm{(i)}}\ $\rho_{ab}=f(\textbf{a};\rho_{ba})$ and $\rho_{ba}=f(\textbf{b};\rho_{ab})$.
\par
{\rm(ii)}\ If $m_a m_b\leq 1$, then $\rho_{_{ab}}=\rho_{_{ba}}=1$.
\par
{\rm(iii)}\ If $m_am_b> 1$ and $\rho_{a}=\rho_{b}$, then
$\rho_{_{ab}}= \rho_{_{ba}}=\rho_{_{a}}$.
\par
{\rm(iv)}\ If $m_am_b> 1$ and $\rho_{a}<\rho_{b}$, then
$\rho_a< \rho_{ab}< \rho_{ba}< \rho_b$.
\end{lemma}
\par
\begin{proof}\ It is easy to see that
$$
\rho_{ab}=\alpha_n(\rho_{ab})=f(\emph{\textbf{a}};
f(\emph{\textbf{b}};\alpha_{n-1}(\rho_{ab})))=f(\emph{\textbf{a}};
\beta_{n-1}(f(\emph{\textbf{b}};\rho_{ab})).
$$
Noting that $\beta_{n-1}(f(\emph{\textbf{b}};\rho_{ab}))\rightarrow \rho_{ba}$ yields $\rho_{ab}=f(\emph{\textbf{a}};\rho_{ba})$. Similarly, $\rho_{ba}=f(\emph{\textbf{b}};\rho_{ab})$. (i) is proved. (ii) follows from $\alpha'(1)=\beta'(1)=m_am_b\leq 1$.
\par
If $m_am_b>1$, then $\alpha'(1)=\beta'(1)=m_am_b>1$ and hence $\rho_{ab}, \rho_{ba}<1$. If further $\rho_a=\rho_b$, then $\rho_a=\rho_b<1$ and
\begin{eqnarray*}
 \alpha(\rho_a)=f(\emph{\textbf{a}};f(\emph{\textbf{b}};\rho_b))
 =f(\emph{\textbf{a}};\rho_b)=f(\emph{\textbf{a}};\rho_a)=\rho_a.
\end{eqnarray*}
Similarly, $\beta(\rho_a)=\rho_a$.
Hence, $\rho_{ab}= \rho_{ba}=\rho_a$. (iii) is proved.
\par
Now, we prove (iv). it is obvious that $\rho_{ab},\ \rho_{ba}<1$ since $m_am_b> 1$. If $\rho_a<\rho_b<1$, by the property of $f(\emph{\textbf{a}};s)$ and $f(\emph{\textbf{b}};s)$,
\begin{eqnarray*}
 \alpha(\rho_a)
  =f(\emph{\textbf{a}};f(\emph{\textbf{b}};\rho_a))
  >f(\emph{\textbf{a}};\rho_a)=\rho_a
\end{eqnarray*}
and
\begin{eqnarray*}
  \alpha(\rho_b)
  =f(\emph{\textbf{a}};f(\emph{\textbf{b}};\rho_b))
  =f(\emph{\textbf{a}};\rho_b)
  <\rho_b.
\end{eqnarray*}
Hence, $\rho_{a}<\rho_{ab}
<\rho_{b}$.
\par
If $\rho_a<\rho_b=1$, by the property of $f(\emph{\textbf{a}};s)$ and $f(\emph{\textbf{b}};s)$,
\begin{eqnarray*}
  \alpha(\rho_a)
  =f(\emph{\textbf{a}};f(\emph{\textbf{b}};\rho_a))
  >f(\emph{\textbf{a}};\rho_a)
  =\rho_a
\end{eqnarray*}
and hence, $\rho_a<\rho_{ab}
<\rho_{b}=1$.
The proof is complete.
\hfill $\Box$
\end{proof}
\par
Similar as in the ordinary Galton-Watson case, we give the following definition.
\par
\begin{definition}\label{def2.b}
An $\emph{\textbf{a}}\!*\!\emph{\textbf{b}}$-Galton-Watson system $\{Z_n: n\geq 0\}$ is called critical, supercritical or subcritical if $m_am_b=1$, $m_am_b>1$ or $m_am_b<1$ respectively.
\end{definition}
\par
Let $\{Z_n: n\geq 0\}$ be an $\emph{\textbf{a}}\!*\!\emph{\textbf{b}}$-Galton-Watson system, where $\emph{\textbf{a}}=\{a_n\}_{n=0}^{\infty}, \emph{\textbf{b}}=\{b_n\}_{n=0}^{\infty}\in \mathscr{P}$.
It is easy to see that $Z_n$ can be rewritten as
\begin{eqnarray}\label{eq2-1}
Z_{n+1}=\begin{cases}
\sum\limits_{i=1}^{Z_n}
\xi_{n,i}^{(a)},\ \mbox{if}\ n\ \mbox{is\ even},\\
\sum\limits_{i=1}^{Z_n}
\xi_{n,i}^{(b)},\ \mbox{if}\ n\ \mbox{is\ odd},\\
\end{cases}
\end{eqnarray}
where $\{\xi_{n,i}^{(a)}:k\geq 0, i\geq 1\}$ and $\{\xi_{n,i}^{(b)}:k\geq 0, i\geq 1\}$ are independent identically distributed random variable sequences with probability distribution $P(\xi_{n,1}^{(a)}=j)=a_j$ and $P(\xi_{n,1}^{(b)}=j)=b_j$ respectively. Moreover, $\{\xi_{n,i}^{(a)}:k\geq 0, i\geq 1\}$ is independent with $\{\xi_{n,i}^{(b)}:k\geq 0, i\geq 1\}$.

\par
Define
\begin{eqnarray}\label{eq2-2}
\begin{cases}
&f_0(\emph{\textbf{ab}};s)=s,\\
&f_{2n+1}(\emph{\textbf{ab}};s)=f_{2n}(\emph{\textbf{ab}};
f(\emph{\textbf{a}};s)),\ \ n\geq 0,\\
&f_{2n+2}(\emph{\textbf{ab}};s)=f_{2n+1}(\emph{\textbf{ab}};
f(\emph{\textbf{b}};s)),\ \ n\geq 0.
\end{cases}
\end{eqnarray}
\par
\begin{theorem}\label{th2.1} Let $\{Z_n: n\geq 0\}$ be an $\textbf{a}\!*\!\textbf{b}$-Galton-Watson system with $Z_0=1$. Then
\par
{\rm{(i)}}\ The probability generating function of $Z_n$ is given by
\begin{eqnarray}\label{eq2-3}
E[s^{Z_n}]=f_n(\textbf{ab};s).
\end{eqnarray}
\par
{\rm{(ii)}}\ The mean and variance of $Z_n$ are given by
\begin{eqnarray}\label{eq2-4}
E[Z_n]=\Gamma_n
\end{eqnarray}
and
\begin{eqnarray}\label{eq2-5}
Var(Z_n)=
\begin{cases}
\frac{\sigma^2m^{k-1}(m^k-1)m_a^2}{(m-1)}+\sigma_a^2m^k,~n=2k+1\\
\frac{\sigma^2m^{k-1}(m^k-1)}{(m-1)},~~n=2k
\end{cases}
\end{eqnarray}
where $\sigma_a^2=\sum\limits_{j=1}^{\infty}j^2a_j-m_a^2$, $\sigma_b^2=\sum\limits_{j=1}^{\infty}j^2b_j-m_b^2$ and $\sigma^2=\sigma_a^2m_b^2+\sigma_b^2m_a$.
\end{theorem}
\par
\begin{proof} We first prove (i). If $n=1$, then
\begin{eqnarray*}
E[s^{Z_1}]=
E[s^{\xi^{(a)}_{0,1}}]=\sum_{k=0}^\infty a_ks^k=f(\emph{\textbf{a}};s)=f_1(\emph{\textbf{ab}};s).
\end{eqnarray*}
\par
If $n=2$, then
\begin{eqnarray*}
E[s^{Z_2}]
=E[s^{\sum\limits_{i=1}^{Z_1}\xi^{(b)}
_{1,i}}
]
=\sum_{j=0}^\infty P(Z_1=j)\cdot E[s^{\sum\limits_{i=1}^j \xi_{1,i}^{(b)}}]
=f(\emph{\textbf{a}};f(\emph{\textbf{b}};s))=f_2(\emph{\textbf{ab}};s).
\end{eqnarray*}
\par
Suppose (\ref{eq2-3}) holds true for $n$. Then, if $n=2k$, we have
\begin{eqnarray*}
E[s^{Z_{2k+1}}]=E[s^{\sum\limits_{i=1}
^{Z_{2k}}
\xi^{(a)}_{2k,i}}]=\sum_{j=0}^\infty P(Z_{2k}=j)
\cdot(E[s^{\xi^{(a)}_{2k,1}}])^j=f_{2k}(\emph{\textbf{ab}};f(\emph{
\textbf{a}};s))
=f_{2k+1}(\emph{\textbf{ab}};s).
\end{eqnarray*}
Similarly, if $n=2k+1$, then $E[s^{Z_{2k+2}}]
=f_{2k+2}(\emph{\textbf{ab}};s)$. Therefore, (i) is proved.
\par
Now we prove (ii). By (\ref{eq2-2}),
\begin{eqnarray*}
E[Z_{2k+1}]=f'_{2k}(\emph{
\textbf{ab}};1)\cdot f'(\emph{\textbf{a}};1)=m_a\cdot f'_{2k}(\emph{\textbf{ab}};1),\ \ k\geq 0,
\end{eqnarray*}
\begin{eqnarray*}
E[Z_{2k+2}]=f'_{2k+1}(\emph{
\textbf{ab}};1)\cdot f'(\emph{\textbf{b}};1)=m_b\cdot f'_{2k+1}(\emph{\textbf{ab}};1),\ \ k\geq 0,
\end{eqnarray*}
which implies (\ref{eq2-4}) since $E[Z_0]=1$ and $E[Z_1]=m_a$. On the other hand,
\begin{eqnarray*}
\frac{d^2f_{_{2k+1}}(\emph{\textbf{ab}};s)}{ds^2}=f''_{2k}(\emph{
\textbf{ab}};f(\emph{\textbf{a}};s))\cdot (f'(\emph{\textbf{a}};s))^2+f'_{2k}(\emph{\textbf{ab}};f(\emph{
\textbf{a}};s))\cdot f''(\emph{\textbf{a}};s),\ \ k\geq 0.
\end{eqnarray*}
Denote $V_{n}=E[Z_n(Z_n-1)]$.
Then
\begin{eqnarray*}
V_{2k+1}=V_{2k}\cdot m_a^2+m^k\cdot (\sigma_a^2-m_a+m_a^2),\ \ k\geq 0.
\end{eqnarray*}
Similarly,
\begin{eqnarray*}
V_{2k+2}=V_{2k+1}\cdot m_b^2+m^km_a\cdot (\sigma_b^2-m_b+m_b^2),\ \ k\geq 0.
\end{eqnarray*}
By the above two equalities,
\begin{eqnarray*}
V_{2k+2}&=& m^2\cdot V_{2k}+m^km^2_b\cdot (\sigma_a^2-m_a+m_a^2)+m^km_a\cdot (\sigma_b^2-m_b+m_b^2)\\
&=& m^2\cdot V_{2k}+m^k(m^2_b\sigma_a^2-mm_b+m^2+m_a\sigma_b^2-m+mm_b)\\
&=& m^2\cdot V_{2k}+m^k(\sigma^2+m^2-m),\ \ k\geq 0,
\end{eqnarray*}
which implies
\begin{eqnarray*}
V_{2k}=\frac{m^{k-1}(m^k-1)\sigma^2}{m-1}+m^{2k}-m^k,\ \ k\geq 0
\end{eqnarray*}
and
\begin{eqnarray*}
V_{2k+1}=\frac{m^{k-1}(m^k-1)m_a^2\sigma^2}{m-1}+m^{2k} m_a^2 -m^km_a+m^k\sigma_a^2,\ \ k\geq 0.
\end{eqnarray*}
Hence,
\begin{eqnarray*}
Var(Z_{2k})=V_{2k}+m^k-m^{2k}
=\frac{m^{k-1}(m^k-1)\sigma^2}{m-1},\  \ k\geq 0
\end{eqnarray*}
and
\begin{eqnarray*}
Var(Z_{2k+1})
=V_{2k+1}+m^km_a-m^{2k}m_a^2
=\frac{m^{k-1}(m^k-1)m_a^2\sigma^2}{m-1}+m^k\sigma_a^2,\  \ k\geq 0.
\end{eqnarray*}
The proof is complete. \hfill $\Box$
\end{proof}
\par
\par
\begin{remark}\label{re2.1}
By Theorem~\ref{th2.1}, it is easy to see that $\{Z_{2n}:n\geq 0\}$ is an ordinary Galton-Watson system in which the generating function of branching mechanism is $f(\emph{\textbf{a}}\emph{\textbf{b}};s)$. Hence, by Athreya~\cite{AK94}, we can obtain the convergence rate of $P(|\frac{Z_{2n+2}}
{Z_{2n}}-m|>\varepsilon\mid Z_0=1)$. However, $\{Z_{2n+1}:n\geq 0\}$ is not a real ordinary Galton-Watson system since $E[s^{Z_3}]=f(\emph{\textbf{a}};f(\emph{\textbf{b}};f(\emph{\textbf{a}};s))$ can not be a composite function of $E[s^{Z_1}]=f(\emph{\textbf{a}};s)$ itself.
\end{remark}
\par
The following lemma is due to Athreya and Ney~\cite{AKN72} and the proof is omitted.
\par
\begin{lemma}\label{le2.3}\ Suppose that $a_0+a_1<1, b_0+b_1<1$.
\par
{\rm{(i)}}\ $\alpha(s)$ and $\beta(s)$ are strictly convex and increasing on $[0,1]$;
\par
{\rm{(ii)}}\ $\alpha_n(s)\uparrow \rho_{ab}$ as $n\to\infty$ for $s\in[0,\rho_{ab})$, while $\alpha_n(s)\downarrow \rho_{ab}$ as $n\to\infty$ for $s\in (\rho_{ab},1)$; Similarly, $\beta_n(s)\uparrow \rho_{ba}$ as $n\to\infty$ for $s\in[0,\rho_{ba})$, while $\beta_n(s)\downarrow \rho_{ba}$ as $n\to\infty$ for $s\in (\rho_{ba},1)$.
\end{lemma}
\par
\begin{theorem}\label{th2.2}
Let $\{Z_n: n\geq 0\}$ be an $\textbf{a}\!*\!\textbf{b}$-Galton-Watson system with $Z_0=1$. Then the extinction probability of $\{Z_n: n\geq 0\}$ is $\rho_{ab}$, which is the smallest nonnegative root of $\alpha(s)=s$.
\end{theorem}
\par
\begin{proof}
By Theorem~\ref{th2.1}, we only need to prove that $\lim\limits_{n\to \infty}f_n(\emph{\textbf{ab}};s)=\rho_{ab}$ for all $s\in (0,1)$. Without loss of generality, we assume $a_0,b_0<1$. Since
\begin{eqnarray*}
&&f_{1}(\emph{\textbf{ab}};s)=f(\emph{\textbf{a}};s)\\
&&f_{2}(\emph{\textbf{ab}};s)
=f(\emph{\textbf{a}};f(\emph{\textbf{b}};s))=\alpha(s).\\
\end{eqnarray*}
Recursively,
\begin{eqnarray*}
&&f_{2n}(\emph{\textbf{ab}};s)=\alpha_n(s),\\
&&f_{2n+1}(\emph{\textbf{ab}};s)=\alpha_n(f(\emph{\textbf{a}};s)).
\end{eqnarray*}
Since $f(\emph{\textbf{a}};s)<1$ for all $s\in [0,1)$, by Lemma~\ref{le2.3}, we know that
$$
\lim_{n\to \infty}f_{2n}(\emph{\textbf{ab}};s)=\lim_{n\to \infty}
f_{2n+1}(\emph{\textbf{ab}};s)=\rho_{ab}.
$$
Hence, $\lim\limits_{n\to \infty}f_n(\emph{\textbf{ab}};s)=\rho_{ab}$ for all $s\in (0,1)$. The proof is complete. \hfill $\Box$
\end{proof}
\par
\vspace{5mm}
 \setcounter{section}{3}
 \setcounter{equation}{0}
 \setcounter{theorem}{0}
 \setcounter{lemma}{0}
 \setcounter{corollary}{0}
 \setcounter{assumption}{0}

\noindent {\large \bf 3. Limit properties}
\vspace{3mm}
\par
In this section, we study the limit properties of $Z_n$. Let $W_n$ be defined in (\ref{eq1-3}), i.e.,
\begin{eqnarray*}
W_n=\Gamma_n^{-1}Z_n,\quad n\geq 0,
\end{eqnarray*}
where $\Gamma_n$ is defined by (\ref{eq1-2}).
\par
The following theorem reveals the properties of $W_n$.
\par
\begin{theorem}\label{th3.1}
Suppose that $m>1,~ \sigma_a^2, \sigma_b^2<\infty$, and $Z_0=1$. Then $W_n$ is an integrable martingale and hence converges to some random variable $W$. Moreover,
\par
{\rm{(i)}}\ $\lim\limits_{n\to\infty}E[(W_n-W)^2]=0$;
\par
{\rm{(ii)}}\ $E[W]=1, Var(W)=\frac{\sigma^2}{m^2-m}$;
\par
{\rm{(iii)}}\ $P(W=0)=\rho_{ab}$.
\end{theorem}
\par
\begin{proof}
By Markov property and (\ref{eq2-1}), we have
\begin{eqnarray*}
&& E[W_{n+1}\mid\sigma(Z_0,Z_1,\cdots,Z_n)]\\
&=& E[\Gamma_{n+1}^{-1}\sum\limits_{i=1}^{Z_n}\xi^{(\emph{\textbf{a}}^{(n)})}_{n,i}\mid Z_n]\\
&=&\Gamma_{n+1}^{-1}Z_nf'(\emph{\textbf{a}}^{(n)};1)\\
&=&\Gamma_n^{-1}Z_n,
\end{eqnarray*}
where $\emph{\textbf{a}}^{(n)}=\emph{\textbf{a}}$ or $\emph{\textbf{b}}$ for $n$ being even or odd respectively, and we have used the fact that $\Gamma_{n+1}^{-1}f'(\emph{\textbf{a}}^{(n)};1)=\Gamma_n^{-1}$ for all $n$. Hence, by Theorem~\ref{th2.1}, $W_n$ is an integrable martingale and hence converges to some random variable $W$.
\par
It follows from (\ref{eq2-5}) that
\begin{equation*}
E[W_n^2]=\Gamma_n^{-2}E[Z_n^2]=
\begin{cases}
\frac{\sigma^2(1-m^{-k})}{m^2-m}+\frac{\sigma_a^2}{m^km_a^2}+1, n=2k+1\\
\frac{\sigma^2(1-m^{-k})}{m^2-m}+1,~n=2k
\end{cases}
\end{equation*}
and hence, $\sup\limits_nE[W_n^2]<\infty$. Now by standard martingale theory (see Doob~\cite{DJ}), (i) and (ii) follow.
\par
If $r=P(W=0)$ then $E[W]=1$ implies $r<1$. Furthermore,
\[
r=\lim_{n\to \infty}P(W_n=0)=\lim_{n\to \infty}P(Z_n=0)=\rho_{ab}.
\]
The proof is complete. \hfill $\Box$
\end{proof}
\par
By the proof of Theorem~\ref{th2.2}, we know that for any $s\in[0,1)$, $f_n(\emph{\textbf{ab}};s)\to \rho_{ab}$ as $n\to\infty$. We shall now consider the convergence rate. Define
\begin{eqnarray}\label{eq3-1}
Q_n(s)=\frac{f_n(\emph{\textbf{ab}};s)-\rho_{ab}}{\gamma_n} ,
\end{eqnarray}
where
\begin{eqnarray*}
\gamma_n=\begin{cases}
[f'(\emph{\textbf{a}};\rho_{ba})\cdot f'(\emph{\textbf{b}};\rho_{ab})]^k,\ & n=2k, \\
[f'(\emph{\textbf{a}};\rho_{ba})\cdot f'(\emph{\textbf{b}};\rho_{ab})]^k\cdot f'(\emph{\textbf{a}},\rho_{ba}),\ & n=2k+1.\\
\end{cases}
\end{eqnarray*}
\par
\begin{theorem}\label{th3.2}
Suppose that $m>1$. Then
\begin{eqnarray}\label{eq3-2}
\lim\limits_{k\to\infty}Q'_{2k}(s)=:Q'(s)\ \ {\rm{and}}\ \lim\limits_{k\to\infty}Q'_{2k+1}(s)=:\tilde{Q}'(s)
\end{eqnarray}
exist for $0\leq s<1$ and $Q'(s)>0$ for all $s\in [0,1)$. Furthermore, $\lim\limits_{s\to \rho_{ab}}Q'(s)=1$ and
\begin{eqnarray}\label{eq3-3}
\tilde{Q}'(s)=Q'(f(\textbf{a};
s))\cdot \frac{f'(\textbf{a};s)}{f'(\textbf{a};\rho_{ba})}.
\end{eqnarray}
\end{theorem}
\par
\begin{proof}
It is easy to see that
\begin{equation}\label{eq3-4}
Q'_n(s)=\frac{f_n'(\emph{\textbf{ab}};s)}{\gamma_n}=\begin{cases}
\frac{\alpha'_k(s)}{[f'(\emph{\textbf{a}};\rho_{ba})\cdot f'(\emph{\textbf{b}};\rho_{ab})]^k},\ & n=2k,\\
\frac{\alpha'_k(f(\emph{\textbf{a}};s))}{[f'(\emph{\textbf{a}};
\rho_{ba})\cdot f'(\emph{\textbf{b}};\rho_{ab})]^k}\cdot \frac{f'(\emph{\textbf{a}};s)}{f'(\emph{\textbf{a}};\rho_{ba})},\ & n=2k+1.\\
\end{cases}
\end{equation}
\par
Since $m=m_am_b=f'(\emph{\textbf{a}};1)\cdot f'(\emph{\textbf{b}};1)=\alpha'(1)>1$, by Theorem 1.11.1 of Arthreya~\cite{AKN72}, we know that
\begin{equation*}
\lim\limits_{k\to\infty}\frac{\alpha'_k(s)}{[f'(\emph{\textbf{a}};
\rho_{ba})\cdot f'(\emph{\textbf{b}};\rho_{ab})]^k}=Q'(s)
\end{equation*}
exists for all $s\in [0,1)$ and $Q'(s)>0$ for all $s\in [0,1)$. Furthermore,
$\lim\limits_{s\to \rho_{ab}}Q'(s)=1$. Hence, $\tilde{Q}(s)$ exists and by (\ref{eq3-4}),
\begin{eqnarray*}
\tilde{Q}'(s)>0,\ s\in [0,1)\ \ {\rm{and}}\ \ \lim_{s\to \rho_{ba}}\tilde{Q}'(s)=1.
\end{eqnarray*}
The proof is complete. \hfill $\Box$
\end{proof}
\par
Now define
\begin{equation*}
Q(s)=\int_{\rho_{_{ab}}}^sQ'(t)dt, ~~for~~0\leq s<1,
\end{equation*}
\begin{equation*}
\tilde{Q}(s)=\int_{\rho_{_{ba}}}^s\tilde{Q}'(t)dt, ~~for~~0\leq s<1.
\end{equation*}
Then we have
\begin{corollary}\label{cor3.1}\ If $m>1$. Then,
\begin{equation*}
\lim\limits_{n\to\infty}Q_{2n}(s)=Q(s)\ \ {\rm{and}}\ \
\lim\limits_{n\to\infty}Q_{2n+1}(s)=\tilde{Q}(s).
\end{equation*}
\end{corollary}
\par
\begin{proof} By the bounded convergence theorem,
\begin{equation*}
Q_{2n}(s)=
Q_{2n}(s)-Q_{2n}(\rho_{ab})=\int_{\rho_{_{ab}}}^sQ'_{2n}(t)dt\to Q(s).
\end{equation*}
and
\begin{equation*}
Q_{2n+1}(s)=
Q_{2n+1}(s)-Q_{2n+1}(\rho_{ba})=\int_{\rho_{_{ba}}}^sQ'_{2n+1}(t)dt
\to \tilde{Q}(s).
\end{equation*}
The proof is complete. \hfill $\Box$
\end{proof}
\par
\begin{theorem}\label{th3.3}
$(Q(s),\tilde{Q}(s))$ is the unique solution of the functional equations
\begin{equation}\label{eq3-5}
\begin{cases}
Q(f(\textbf{a};s))=f'(\textbf{a};\rho_{ba})\cdot \tilde{Q}(s),  \\
\tilde{Q}(f(\textbf{b};s))=f'(\textbf{b};\rho_{ab})\cdot Q(s),  \\
\end{cases}
\end{equation}
subject to
\begin{equation}\label{eq3-6}
Q(\rho_{ab})=0,~\tilde{Q}(\rho_{ba})=0 ~and~\lim\limits_{s\to \rho_{ab}}Q'(s)=1, ~\lim\limits_{s\to \rho_{ba}}\tilde{Q}'(s)=1.
\end{equation}
\end{theorem}
\par
 \begin{proof}\ Substituting $f(\emph{\textbf{a}};s)$ or $f(\emph{\textbf{b}};s)$ for $s$ in the definition of $Q_n(s)$,
\begin{equation*}
Q_{2k}(s)=\frac{f_{2k}(\emph{\textbf{ab}};s)-\rho_{ab}}{{\gamma}_{2k}}
=\frac{f_{2k-1}(\emph{\textbf{ab}};f(\emph{\textbf{b}};s))-\rho_{ab}}
{{\gamma}_{2k-1}\cdot f'(\emph{\textbf{b}};\rho_{ab})}=\frac{Q_{2k-1}(f(\emph{\textbf{b}};
s))}{
f'(\emph{\textbf{b}};\rho_{ab})}.
\end{equation*}
Taking limits on both sides yields the second equality of (\ref{eq3-5}). The first equality of (\ref{eq3-5}) can be similarly proved.
\par
 As for the uniqueness, note that if $(Q(s),\tilde{Q}(s))$ and $(H(s),\tilde{H}(s))$ are two solutions of (\ref{eq3-5}) subject to (\ref{eq3-6}), then
\begin{eqnarray}\label{eq3-7}
|\tilde{Q}(s)-\tilde{H}(s)|
&=& f'(\emph{\textbf{a}};\rho_{ba})^{-1}|Q(f(\emph{\textbf{a}};s))
-H(f(\emph{\textbf{a}};s))|\nonumber\\
&=& f'(\emph{\textbf{a}};\rho_{ba})^{-1}f'(\emph{\textbf{b}};\rho_{ab})
^{-1}
\cdot|\tilde{Q}(f(\emph{\textbf{b}};f(\emph{\textbf{a}};s))) -\tilde{H}(f(\emph{\textbf{b}};f(\emph{\textbf{a}};s)))|\nonumber\\
&=&{\gamma}_{_{2k}}^{-1}|\tilde{Q}(f_{_{2k}}(\emph{\textbf{ba}};s))
-\tilde{H}(f_{_{2k}}(\emph{\textbf{ba}};s))|\nonumber\\
&\leq & |\hat{Q}_{2k}(s)|\cdot \left\{|1-\frac{\tilde{Q}(f_{2k}(\emph{\textbf{ba}};s))}
{f_{2k}(\emph{\textbf{ba}};s)-\rho_{ba}}|
+|1-\frac{\tilde{H}(f_{2k}(\emph{\textbf{ba}};s))}
{f_{2k}(\emph{\textbf{ba}};s)-\rho_{ba}}|\right\},
\end{eqnarray}
where
\begin{equation*}
\hat{Q}_n(s)=\frac{f_n(\emph{\textbf{ba}};s)-\rho_{ba}}{\gamma'_n} ,
\end{equation*}
\begin{eqnarray*}
\gamma'_n=\begin{cases}
[f'(\emph{\textbf{b}};\rho_{ab})\cdot f'(\emph{\textbf{a}};\rho_{ba})]^k, & n=2k, \\
[f'(\emph{\textbf{b}};\rho_{ab})\cdot f'(\emph{\textbf{a}};\rho_{ba})]^kf'(\emph{\textbf{b}};\rho_{ab}), & n=2k+1. \\
\end{cases}
\end{eqnarray*}
\par
Now for any $s\in [0,1)$, $f_{2k}(\emph{\textbf{ba}};s)\to \rho_{ba}$ and
\begin{equation*}
\lim\limits_{k\to \infty}\frac{\tilde{Q}(f_{2k}(\emph{\textbf{ba}};s))}{f_{2k}(\emph{
\textbf{ba}};s)
-\rho_{ba}}=\lim\limits_{s\to \rho_{ba}}\tilde{Q}'(s)=1,\ \ \ \lim\limits_{k\to \infty}\frac{\tilde{H}(f_{2k}(\emph{\textbf{ba}};s))}{f_{2k}(\emph{
\textbf{ba}};s)
-\rho_{ba}}=\lim\limits_{s\to \rho_{ba}}\tilde{H}'(s)=1.
\end{equation*}
Similarly, we have
\begin{equation}\label{eq3-8}
|Q(s)-H(s)|\leq|\hat{Q}_{2k-1}(s)|\cdot \left\{|1-\frac{Q(f_{2k-1}(\emph{\textbf{ab}};s))}
{f_{2k-1}(\emph{\textbf{ab}};s)-\rho_{ab}}|
+|1-\frac{H(f_{2k-1}(\emph{\textbf{ab}};s))}{f_{2k-1}(\emph{
\textbf{ab}};s)
-\rho_{ab}}|\right\},
\end{equation}
and
\begin{equation*}
\lim\limits_{k\to \infty}\frac{Q(f_{2k-1}(\emph{\textbf{ab}};s))}
{f_{2k-1}(\emph{\textbf{ab}};s)-\rho_{ab}}=\lim\limits_{s\to \rho_{ab}}Q'(s)=1,\ \ \ \lim\limits_{k\to \infty}\frac{H(f_{2k-1}(\emph{\textbf{ab}};s))}
{f_{2k-1}(\emph{\textbf{ab}};s)-\rho_{ab}}=\lim\limits_{s\to \rho_{ab}}H'(s)=1.
\end{equation*}
On the other hand, by a similar argument of Theorem~\ref{th3.2} and Corollary~\ref{cor3.1}, we know that $\lim_{k\to \infty}\hat{Q}_{2k}(s)$ and $\lim_{k\to \infty}\hat{Q}_{2k+1}(s)$ exist and are finite for all $s\in [0,1)$. Therefore, by (\ref{eq3-7}) and (\ref{eq3-8}), we have $Q(s)=H(s)$ and $\tilde{Q}(s)=\tilde{H}(s)$ for all $s\in [0,1)$. The proof is complete.
\hfill $\Box$
\end{proof}
\par
Since $Q(s)$ and $\tilde{Q}(s)$ are limits of power series, we may rewrite them as
\begin{eqnarray}\label{eq3-9}
Q(s)=\sum_{k=0}^{\infty}q_k s^k,\ \ s\in [0,1)\ \ \ {\rm{and}}\ \ \ \tilde{Q}(s)=\sum_{k=0}^{\infty}\tilde{q}_k s^k,\ \ s\in [0,1).
\end{eqnarray}
\par
\vspace{5mm}
 \setcounter{section}{4}
 \setcounter{equation}{0}
 \setcounter{theorem}{0}
 \setcounter{lemma}{0}
 \setcounter{corollary}{0}
 \setcounter{assumption}{0}

\noindent {\large \bf 4. Large deviation}
\vspace{3mm}
\par
In this section, we will discuss the large deviation rates of $\{Z_n:n\geq 0\}$.
By Theorem~\ref{th2.1}, we know that $E\{s^{Z_n}\mid Z_0=1\}=f_n(\emph{\textbf{ab}};s)$. Therefore, we first study the convergence property of $f_n(\emph{\textbf{ab}};s)$ and its inverse as $n\to\infty$.
\par
From now on, we assume that
\begin{eqnarray*}
a_0=b_0=0, a_j, b_j\neq 1, \ \forall j\geq 0\  {\rm{and}}\ \ m_a=f'(\emph{\textbf{a}};1), m_b=f'(\emph{\textbf{b}};1)<\infty.
\end{eqnarray*}
\par
\begin{proposition}\label{prop3.1} Let $a_0=b_0=0,a_1,b_1>0$. Then $\rho_a=\rho_{ab}=\rho_{ba}=\rho_b=0$ and there exist $0\leq q_j<\infty$, $0\leq \tilde{q}_j<\infty$ $(j\geq 1)$ with $q_1=\tilde{q}_1=1$ such that
\begin{equation}\label{eq4-1}
\lim\limits_{n\to \infty}\frac{f_{2n}(\textbf{ab};s)}{(a_1b_1)^n}=\sum_{j=1}^{\infty}
q_j s^j=Q(s)<\infty, \ \ s\in [0,1)
\end{equation}
and
\begin{equation}\label{eq4-2}
\lim\limits_{n\to \infty}\frac{f_{2n+1}(\textbf{ab};s)}{(a_1b_1)^na_1}
=\sum_{j=1}^{\infty}\tilde{q}_j s^j=\tilde{Q}(s)<\infty, \ \ s\in [0,1).
\end{equation}
Furthermore, $(Q(s),\tilde{Q}(s))$ is the unique solution of the functional equations
\begin{equation*}
\begin{cases}
Q(f(\textbf{a};s))=a_1\tilde{Q}(s)  \\
\tilde{Q}(f(\textbf{b};s))=b_1 Q(s)   \\
\end{cases}
\end{equation*}
subject to
\begin{equation*}
Q(0)=0,~\tilde{Q}(0)=0, Q(1),\tilde{Q}(1)=\infty\ \ {\rm{and}}\ \ Q(s), \tilde{Q}(s)<\infty,\ s\in (0,1).
\end{equation*}
Consequently, for all $1\leq i,j<\infty$,
\begin{equation}\label{eq4-3}
\lim\limits_{n\to \infty}\frac{P(Z_{2n}=j|Z_0=i)}{(a_1b_1)^{in}}=q^{*(i)}_j
\end{equation}
and
\begin{equation}\label{eq4-4}
\lim\limits_{n\to \infty}\frac{P(Z_{2n+1}=j|Z_0=i)}{(a_1b_1)^{in}a_1^i}=\tilde{q}
^{*(i)}_j,
\end{equation}
where $q^{*(i)}_j$, $\tilde{q}^{*(i)}_j$ satisfy $\sum_{j=1}^{\infty}q^{*(i)}_js^j=Q^i(s)$ and $\sum_{j=1}^{\infty}\tilde{q}^{*(i)}_js^j=\tilde{Q}^i(s)$ for $s\in [0,1)$.
\end{proposition}
\par
\begin{proof}\ All the assertions excepting (\ref{eq4-3}) and (\ref{eq4-4}) follow directly from Theorems~\ref{th3.2}-\ref{th3.3} and Corollary~\ref{cor3.1}. while (\ref{eq4-3}) and (\ref{eq4-4}) follow from (\ref{eq4-1}) and (\ref{eq4-2}) since $E[s^{Z_n}|Z_0=i]=(E[s^{Z_n}|Z_0=1])^i=f^i_n(\emph{\textbf{ab}};s)$.
\hfill $\Box$
\end{proof}
\par
Now, denote $r_a:=\max\{s:f(\emph{\textbf{a}};s)<\infty\}$ and $r_b:=\max\{s:f(\emph{\textbf{b}};s)<\infty\}$. Obviously, $r_a,r_b\geq 1$. Let $g(\emph{\textbf{a}};s)$ and $g(\emph{\textbf{b}};s)$ be the inverse functions of $f(\emph{\textbf{a}};s)$ and of $f(\emph{\textbf{b}};s)$ respectively, which are defined by
\begin{equation*}
f(\Delta;g(\Delta;s))=s ~~\text{for}~~0\leq s<\infty,~~\Delta=\emph{\textbf{a}},\emph{\textbf{b}}.
\end{equation*}
\par
It can be easily seen that $g(\emph{\textbf{a}};s)$ and $g(\emph{\textbf{b}};s)$ are well defined on $[0,f(\emph{\textbf{a}};r_a)]$ and $[0,f(\emph{\textbf{b}};r_b)]$ respectively. Moreover, $g(\emph{\textbf{a}};s),g(\emph{\textbf{b}};s)\geq s$ for $s\in [0,1]$ and $g(\emph{\textbf{a}};s)\leq s$ for $s\in [1,r_a]$ (respectively, $g(\emph{\textbf{b}};s)\leq s$ for $s\in [1,r_b]$). Let $g_n(\emph{\textbf{ab}};s)$ and $g_n(\emph{\textbf{ba}};s)$ be the inverse functions of $f_n(\emph{\textbf{ab}};s)$ and $f_n(\emph{\textbf{ba}};s)$ respectively. It is easy to see that
\begin{eqnarray*}
\begin{cases}
g_{1}(\emph{\textbf{ab}};s)=g(\emph{\textbf{a}};s)\\
g_{2n}(\emph{\textbf{ab}};s)=g(\emph{\textbf{b}};g_{2n-1}(\emph{
\textbf{ab}};
s)),\ \ n\geq 1\\
g_{2n+1}(\emph{\textbf{ab}};s)=g(\emph{\textbf{a}};g_{2n}(\emph{
\textbf{ab}};
s)),\ \ n\geq 0\\
\end{cases}
\quad and\ \quad
\begin{cases}
g_{1}(\emph{\textbf{ba}};s)=g(\emph{\textbf{b}};s)\\
g_{2n}(\emph{\textbf{ba}};s)=g(\emph{\textbf{a}};
g_{2n-1}(\emph{\textbf{ba}};
s)),\ \ n\geq 1\\
g_{2n+1}(\emph{\textbf{ba}};s)=g(\emph{\textbf{b}};
g_{2n}(\emph{\textbf{ba}};
s)),\ \ n\geq 0.\\
\end{cases}
\end{eqnarray*}
Moreover, $g_n(\emph{\textbf{ab}};s)$ and $g_n(\emph{\textbf{ba}};s)$ are nondecreasing with $n$ for $s\in [0,1]$ and nonincreasing with $n$ for $s\in [1,f(\emph{\textbf{a}};s_0)]$ and $s\in [1,f(\emph{\textbf{b}};s_0)]$ respectively.
\par
 The next proposition shows that the rate of convergence of $g_n(\emph{\textbf{ab}};\cdot)$ is geometric.
\par
\begin{proposition}\label{prop4.2}
 Let $f(\textbf{a};s_0)<\infty$ for some $s_0>1$. Then, for $1\leq s\leq f(\textbf{a};s_0)$, we have $g_n(\textbf{ab};s),~g_n(\textbf{ba};s)\downarrow 1$ and
\begin{equation}\label{eq4-5}
R_n(s):= \Gamma_n\cdot(g_n(\textbf{ab};s)-1)\downarrow R(s),
\end{equation}
\begin{equation}\label{eq4-6}
\tilde{R}_n(s):= \tilde{\Gamma}_n\cdot(g_n(\textbf{ba};s)-1)\downarrow \tilde{R}(s),
\end{equation}
where
\begin{eqnarray}\label{eq4-7}
\tilde{\Gamma}_n=\begin{cases}m^{\frac{n}{2}},\ & n\ \rm{is\ even},\\
m^{\frac{n-1}{2}}m_b,\ & n\ \rm{is\ odd}.
\end{cases}
\end{eqnarray}
Moreover, $(R(\cdot),\tilde{R}(\cdot))$ is the unique solution of the functional equations
\begin{equation}\label{eq4-8}
\begin{cases}
R(f(\textbf{a};s))=m_a\tilde{R}(s),\\
\tilde{R}(f(\textbf{b};s))=m_bR(s),\\
\end{cases}\ \ s\in [1,f(\textbf{a};s_0)]
\end{equation}
subject to
\begin{eqnarray}\label{eq4-9}
\begin{cases}
&0<R(s),~\tilde{R}(s)<\infty\ ~for\ ~s\in [1,f(\textbf{a};s_0)],\\
&R(1)=\tilde{R}(1)=0,\ \ R'(1)=\tilde{R}'(1)=1.
\end{cases}
\end{eqnarray}
\end{proposition}
The proof of Proposition~\ref{prop4.2} is similar to that of Theorem~\ref{th3.3} and is omitted.
\par
We are now at the position to discuss the convergence rates in (\ref{eq1-4})-(\ref{eq1-6}). For convenience, we will assume $Z_0=1$ from now on. The following theorem presents the convergence rates in (\ref{eq1-4}).
\par
\begin{theorem}\ \label{th4.1}
If $a_0=b_0=0, a_1,b_1>0$ and $f(\textbf{a};s_0)+f(\textbf{b};s_0)<\infty$ for some $s_0>1$. Let $\varepsilon>0$. Then there exists $\lambda\in(0,1)$ such that
\begin{equation}\label{eq4-10}
\phi_a(n,\varepsilon)=P(|\overline{X}_n-m_a|>\varepsilon) \\
=o(\lambda^n)~~as~~n\to\infty,
\end{equation}
\begin{equation}\label{eq4-11}
\phi_b(n,\varepsilon)=P(|\overline{Y}_n-m_b|> \varepsilon)=o(\lambda^n)~~as~~n\to\infty,\\
\end{equation}
where $\overline{X}_n=\frac{\sum_{i=1}^nX_i}{n}$ being the mean of $n$ i.i.d. r.v.s $\{X_i\}$ with distribution $\{a_j\}$, and $\overline{Y}_n=\frac{\sum_{i=1}^nY_i}{n}$ being the mean of i.i.d. r.v.s $\{Y_i\}$ with distribution $\{b_j\}$. Furthermore,
\begin{equation}\label{eq4-12}
\lim\limits_{n\to\infty}\frac{1}{a_1^nb_1^n}
P(|\frac{Z_{2n+1}}{Z_{2n}}-m_a|>\varepsilon)=\sum_{j=1}^{\infty}\phi_a(j,\varepsilon)q_j<\infty
\end{equation}
and
\begin{equation}\label{eq4-13}
\lim\limits_{n\to\infty}\frac{1}{a_1^nb_1^{n-1}}P(|\frac{Z_{2n}}
{Z_{2n-1}}-m_b|>\varepsilon)=\sum_{j=1}^{\infty}\phi_b(j,\varepsilon)\tilde{q}_j<\infty,
\end{equation}
where
$\{q_j\}$ and $\{\tilde{q}_j\}$ are defined via $Q(s)=\sum_{j=0}^{\infty}q_js^j$ and $\tilde{Q}(s)=\sum_{j=0}^{\infty}\tilde{q}_j s^j\ ~ (0\leq s<1)$, being the unique solution of functional equations
\begin{eqnarray*}
\begin{cases}
&Q(f(\textbf{a};s))=a_1 \tilde{Q}(s)  \\
&\tilde{Q}(f(\textbf{b};s))=b_1 Q(s)   \\
\end{cases}
\end{eqnarray*}
subject to
\begin{equation*}
Q(0)=0,~\tilde{Q}(0)=0 ~and~\lim\limits_{s\to 0}Q'(s)=1,~\lim\limits_{s\to 0}\tilde{Q}'(s)=1.
\end{equation*}
\end{theorem}
\par
\begin{proof}\ First note that
\begin{equation*}
\begin{split}
\phi_a(n,\varepsilon)=
&\ P(|\overline{X}_n-m_a|>\varepsilon)\\
\leq &\ P(\alpha^{\sum_{i=1}^nX_i}>\alpha^{n(m_a+\varepsilon)})
+P(\beta^{\sum_{i=1}^nX_i}>\beta^{n(m_a-\varepsilon)})\\
\leq &\ [\alpha^{-(m_a+\varepsilon)}f(\emph{\textbf{a}};\alpha)]^n
+[\beta^{-(m_a-\varepsilon)}f(\emph{\textbf{a}};\beta)]^n
\end{split}
\end{equation*}
for any $\alpha>1$ and $\beta<1$. To prove (\ref{eq4-10}),
we only need to show that $\alpha_0^{-(m_a+\varepsilon)}f(\emph{\textbf{a}};\alpha_0)<1$, $\beta_0^{-(m_a-\varepsilon)}f(\emph{\textbf{a}};\beta_0)<1$ for some $\alpha_0>1$ and $\beta_0<1$. Indeed, consider
\begin{equation*}
F(\alpha)=f(\emph{\textbf{a}};\alpha)-\alpha^{m_a+\varepsilon}.
\end{equation*}
It is easy to know $F(1)=0$, and $F'(\alpha)=f'(\emph{\textbf{a}};\alpha)
-\alpha^{m_a+\varepsilon-1}(m_a+\varepsilon)$. When $\alpha\downarrow1$, we have $F'(\alpha)\to m_a-(m_a+\varepsilon)<0$.
Thus, there exists $\alpha_0>1$ such that $f(\emph{\textbf{a}};\alpha_0)<\alpha_0^{m_a+\varepsilon}$, and hence $\alpha_0^{-(m_a+\varepsilon)}f(\emph{\textbf{a}};\alpha_0)<1$.
Similarly, consider
\begin{equation*}
G(\beta)=f(\emph{\textbf{a}};\beta)-\beta^{m_a-\varepsilon}.
\end{equation*}
Then we have $G(1)=0$, and $G'(\beta)=f'(\emph{\textbf{a}};\beta)-\beta^{m_a-\varepsilon-1}
(m_a-\varepsilon)$. When $\beta\uparrow1$, we have $G'(\beta)\to m_a-(m_a-\varepsilon)>0$. Thus, there exists $\beta_0<1$ such that $\beta_0^{-(m_a-\varepsilon)}f(\emph{\textbf{a}};\beta_0)<1$.
(\ref{eq4-10}) is proved. Similarly, (\ref{eq4-11}) holds true.
\par
Next prove (\ref{eq4-12}) and (\ref{eq4-13}). By branching property,
\begin{equation*}
\begin{split}
P(|\frac{Z_{2n+1}}{Z_{2n}}-m_a|> \varepsilon)&=\sum_{j=1}^{\infty}P(Z_{2n}=j)
P(|\frac{Z_{2n+1}}{Z_{2n}}-m_a|>\varepsilon \mid Z_{2n}=j)\\
&=\sum_{j=1}^{\infty}P(Z_{2n}=j)P(|\frac{\sum_{i=1}^jX_i}{j}-m_a|> \varepsilon)\\
&=\sum_{j=1}^{\infty}P(Z_{2n}=j)\phi_a(j,\varepsilon).
\end{split}
\end{equation*}
Let
$h_{2n}(j)=\phi_a(j,\varepsilon)P(Z_{2n}=j) a_1^{-n}b_1^{-n}$.
Take $C$ such that $\phi_a(j,\varepsilon)\leq C\lambda^j$. Then
\begin{equation*}
h_{2n}(j)\leq C\lambda^jP(Z_{2n}=j)a_1^{-n}b_1^{-n}=:r_{2n}(j).
\end{equation*}
Hence
\begin{equation*}
\sum_{j=1}^{\infty} r_{2n}(j)=\sum_{j=1}^{\infty} C\lambda^jP(Z_{2n}=j)a_1^{-n}b_1^{-n}
=Cf_{2n}(\emph{\textbf{ab}};\lambda)a_1^{-n}b_1^{-n}
\to CQ(\lambda)< \infty.
\end{equation*}
Thus $\sum_{j=1}^{\infty} h_{2n}(j)<\infty$ and $h_{2n}(j)\to \phi_a(j,\varepsilon)q_j$. Secondly,
\begin{equation*}
\begin{split}
P(|\frac{Z_{2n}}{Z_{2n-1}}-m_b|>\varepsilon)&=\sum_{j=1}^{\infty}P(Z_{2n-1}=j)
P(|\frac{Z_{2n}}{Z_{2n-1}}-m_b|>\varepsilon\mid Z_{2n-1}=j)\\
&=\sum_{j=1}^{\infty}P(Z_{2n-1}=j)P(|\frac{\sum_{i=1}^jY_i}{j}-m_b|> \varepsilon)\\
&=\sum_{j=1}^{\infty}P(Z_{2n-1}=j)\phi_b(j,\varepsilon).
\end{split}
\end{equation*}
Let $h_{2n-1}(j)=\phi_b(j,\varepsilon)P(Z_{2n-1}=j)a_1^{-n}b_1^{-n+1}
$. Take $C$ such that $\phi_b(j,\varepsilon)\leq C\lambda^j$. Then
\begin{equation*}
h_{2n-1}(j)\leq C\lambda^jP(Z_{2n-1}=j)a_1^{-n}b_1^{-n+1}=:r_{_{2n-1}}(j).
\end{equation*}
Hence
\begin{equation*}
\sum_{j=1}^{\infty}r_{2n-1}(j)=\sum_{j=1}^{\infty} C\lambda^jP(Z_{2n-1}=j)a_1^{-n}b_1^{-n+1}
=Cf_{2n-1}(ab,\lambda)a_1^{-n}b_1^{-n+1}
\to C\tilde{Q}(\lambda)< \infty.
\end{equation*}
Thus $\sum_{j=1}^{\infty} h_{2n-1}(j)<\infty$ and $h_{2n-1}(j)\to \phi_b(j,\varepsilon)\tilde{q}_j$. The proof is complete. \hfill $\Box$
\end{proof}
\par
The next theorem and corollary establish (\ref{eq4-12}),(\ref{eq4-13}) under conditions weaker than $f(\emph{\textbf{a}};s_0)+f(\emph{\textbf{b}};s_0)<\infty$ for some $s_0>1$.
\par
\begin{theorem}\label{th4.2}\ Assume $a_0=b_0=0, a_1, b_1>0$ and that there exist constants $C_\varepsilon$ and $r>0$ such that $a_1b_1(m_am_b)^r>1$, $\phi_a(k,\varepsilon),~\phi_b(k,\varepsilon)\leq C_\varepsilon/k^r$ for all k, where $\phi_a(k,\varepsilon),~\phi_b(k,\varepsilon)$ are defined in $(\ref{eq4-10}), (\ref{eq4-11})$. Then $(\ref{eq4-12}), (\ref{eq4-13})$ hold.
\end{theorem}
\par
\begin{proof}\ Notice that
\begin{equation*}
P(|\frac{Z_{2n+1}}{Z_{2n}}-m_a|>\varepsilon)
=\sum\limits_{k=1}^{\infty}\phi_a(k,\varepsilon)P(Z_{2n}=k).
\end{equation*}
By assumption,
\begin{equation*}
h_{2n}(k):=\frac{\phi_a(k,\varepsilon)P(Z_{2n}=k)}{a_1^nb_1^n}
\leq\frac{C_\varepsilon}{k^r}\frac{P(Z_{2n}=k)}{a_1^nb_1^n}=: h'_{2n}(k),~~say.
\end{equation*}
By (\ref{eq4-1}),
\begin{equation*}
\begin{split}&
h_{2n}(k)\to q_k\phi_a(k,\varepsilon)=: h(k),~~say\\&
h_{2n}'(k)\to C_\varepsilon \frac{q_k}{k_r}.
\end{split}
\end{equation*}
If we show that
\[
\sum\limits_{k=1}^{\infty}h'_{2n}(k)\to \sum\limits_{k=1}^{\infty}C_\varepsilon\frac{q_k}{k^r}<\infty,
\]
then by a slight modification of the Lebesque dominated convergence theorem, we get that
\[
(a_1^{-n}b_1^{-n})\cdot P(|\frac{Z_{2n+1}}{Z_{2n}}-m_a|>\varepsilon)
=\sum\limits_{k=1}^{\infty}h_{2n}(k)\to\sum\limits_{k=1}^{\infty}
h(k)<\infty.
\]
However,
\[
\sum\limits_{k=1}^{\infty}\frac{1}{k^r}\frac{P(Z_{2n}=k)}{a_1^nb_1^n}
=\frac{E[Z_{2n}^{-r}]}{a_1^nb_1^n}.
\]
For any nonnegative r.v. $X$ and $0<p<\infty$,
\begin{equation*}
\begin{split}
E[X^{-p}]
=E[\frac{1}{\Gamma(p)}\int_0^\infty e^{-tX}t^{p-1}dt]
=\frac{1}{\Gamma(p)}\int_0^\infty E[e^{-tX}]t^{p-1}dt.
\end{split}
\end{equation*}
Therefore,
\begin{equation*}
\begin{split}
\frac{E[Z_{2n}^{-r}]}{a_1^nb_1^n}
=\frac{1}{\Gamma(r)}\int_0^\infty \frac{f_{2n}(\emph{\textbf{ab}};e^{-t})}{a_1^nb_1^n}t^{r-1}dt
=\frac{1}{\Gamma(r)}\int_0^1 \frac{f_{2n}(\emph{\textbf{ab}};s)}{a_1^nb_1^n}k(s)ds,
\end{split}
\end{equation*}
where $\Gamma(\cdot)$ is the Gamma-function and
\begin{eqnarray*}
k(s)=\frac{|\log s|^{r-1}}{s}.
\end{eqnarray*}
Since $f_{2n}(\emph{\textbf{ab}};s)/(a_1^nb_1^n)\uparrow Q(s)$, by the monotone convergence theorem
\[
\Gamma(r)\frac{E[Z_{2n}^{-r}]}{a_1^nb_1^n}\uparrow\int_0^1 Q(s) k(s)ds.
\]
So the proof of (\ref{eq4-12}) will be complete if we show $\int_0^1Q(s)k(s)ds<\infty$. Denote $l(s):=g(\emph{\textbf{b}};g(\emph{\textbf{a}};s))$. Then $l(s):=g(\emph{\textbf{b}};g(\emph{\textbf{a}};s))$ is the inverse of $\hat{l}(s):=f(\emph{\textbf{a}};f(\emph{\textbf{b}};s))$. Let $l_m(s)$ and $\hat{l}_m(s)$ be the $m$'th iterate of $l(s)$ and $\hat{l}(s)$, respectively. Then, $l_m(\hat{l}_m(s))=s,~l_{m+1}(s)\geq l_m(s)$ and for $0<s<1$, $l_m(s)\uparrow 1$ and $\hat{l}_m(s)\downarrow 0$. Fix $0<t_0<1$. Then $t_m=l_m(t_0)\uparrow 1$. Also since $Q(s)$ satisfies(\ref{eq3-5}), (\ref{eq3-6}),
\begin{equation*}
\begin{split}
I_m&
=\int_{t_m}^{t_{m+1}}Q(s)k(s)ds=\int_{t_m}^{t_{m+1}}
\frac{\tilde{Q}(f(\emph{\textbf{b}};s))}{b_1}k(s)ds\\&
=\int_{f(\emph{\textbf{b}};t_m)}^{f(\emph{\textbf{b}};
t_{m+1})}\tilde{Q}(u)
\frac{k(g(\emph{\textbf{b}};u))g'(\emph{\textbf{b}};u)du}
{b_1}=\int_{t_{m-1}}^{t_m}Q(s)
\frac{k(l(s))l'(s)}{a_1b_1}ds\\&
=\int_{t_{m-1}}^{t_m}Q(s)k(s)\frac{k(l(s))l'(s)}{a_1b_1k(s)}ds.
\end{split}
\end{equation*}
Since $l'(s)=1/\hat{l}'(s)$ and $|\log s|/(1-s)\to 1$ as $s\uparrow1$,
$$
\frac{k(l(s))l'(s)}{a_1b_1k(s)}\to \frac{1}{a_1b_1m_a^rm_b^r},
$$
where $m_a=f'(\emph{\textbf{a}};1)$, $m_b=f'(\emph{\textbf{b}};1)$. Thus if $a_1b_1m_a^rm_b^r>1$, then for any $0<(a_1b_1m_a^rm_b^r)^{-1}<\lambda<1$, there exists an $m_0$ such that $k(l(s))l'(s)/(a_1b_1k(s))<\lambda$ for $s\geq l_{m_0}(t_0)$. Thus, $I_m\leq\lambda I_{m-1}$ for $m\geq m_0+2$. Hence,
\begin{equation*}
\begin{split}
\sum\limits_{m=m_0+2}^\infty I_m\leq I_{m_0+1}\sum\limits_{j=1}^\infty\lambda^j<\infty.
\end{split}
\end{equation*}
Therefore,
$$
\int_0^1Q(s)k(s)ds\leq \int_0^{t_{m_0}}Q(s)k(s)ds+\int_{t_{m_0}}^1Q(s)k(s)ds<\infty.
$$
\par
On the other hand,
\begin{equation*}
P(|\frac{Z_{2n}}{Z_{2n-1}}-m_b|>\varepsilon)
=\sum\limits_{k=1}^{\infty}\phi_b(k,\varepsilon)P(Z_{2n-1}=k).
\end{equation*}
By assumption,
\begin{equation*}
h_{2n-1}(k):=\frac{\phi_b(k,\varepsilon)P(Z_{2n-1}=k)}
{a_1^{n}b_1^{n-1}}\leq\frac{C_\varepsilon}{k^r}
\frac{P(Z_{2n-1}=k)}{a_1^{n}b_1^{n-1}}=: h'_{2n-1}(k),~~say.
\end{equation*}
By (\ref{eq4-2}),
\begin{equation*}
\begin{split}&
h_{2n-1}(k)\to \tilde{q}_k\phi_b(k,\varepsilon)=: h(k),~~say,\\&
h'_{2n-1}(k)\to C_\varepsilon \frac{\tilde{q}_k}{k_r}.
\end{split}
\end{equation*}
If we show that
\begin{eqnarray*}
\sum\limits_kh'_{2n-1}(k)\to \sum\limits_kC_\varepsilon\frac{\tilde{q}_k}{k^r}<\infty,
\end{eqnarray*}
then by a slight modification of the Lebesque dominated convergence theorem, we get that
\begin{eqnarray*}
a_1^{-n}b_1^{-n+1}\cdot P(|\frac{Z_{2n}}{Z_{2n-1}}-m_b|>\varepsilon)
=\sum\limits_{k=1}^{\infty}h_{2n-1}(k)\to
\sum\limits_{k=1}^{\infty}h(k)<\infty.
\end{eqnarray*}
However,
\begin{eqnarray*}
\sum\limits_{k=1}^{\infty}\frac{1}{k^r}
\frac{P(Z_{2n-1}=k)}{a_1^nb_1^{n-1}}=\frac{E[Z_{2n-1}^{-r}]}
{a_1^{n}b_1^{n-1}}.
\end{eqnarray*}
For any nonnegative r.v. $X$ and $0<p<\infty$,
\begin{equation*}
\begin{split}
EX^{-p}
=E[\frac{1}{\Gamma(p)}\int_0^\infty e^{-tX}t^{p-1}dt]
=\frac{1}{\Gamma(p)}\int_0^\infty E[e^{-tX}]t^{p-1}dt.
\end{split}
\end{equation*}
Therefore,
\begin{equation*}
\begin{split}
\frac{E[Z_{2n-1}^{-r}]}{a_1^nb_1^{n-1}}
=\frac{1}{\Gamma(r)}\int_0^\infty \frac{f_{2n-1}(\emph{\textbf{ab}};e^{-t})}{a_1^nb_1^{n-1}}t^{r-1}dt
=\frac{1}{\Gamma(r)}\int_0^1 \frac{f_{2n-1}(\emph{\textbf{ab}};s)}{a_1^nb_1^{n-1}}k(s)ds,
\end{split}
\end{equation*}
where
\begin{eqnarray*}
k(s)=\frac{|\log s|^{r-1}}{s}.
\end{eqnarray*}
Since ${f_{2n-1}(\emph{\textbf{ab}};s)}/{a_1^nb_1^{n-1}}\uparrow \tilde{Q}(s)$, by the monotone convergence theorem
\begin{eqnarray*}
\Gamma(r)\frac{E[Z_{2n-1}^{-r}]}{a_1^nb_1^{n-1}}\uparrow\int_0^1 \tilde{Q}(s) k(s)ds,\ \ n\to \infty.
\end{eqnarray*}
So the proof (\ref{eq4-13}) will be complete if we show $\int_0^1\tilde{Q}(s)k(s)ds<\infty$. Denote $\eta(s):=g(\emph{\textbf{a}};g(\emph{\textbf{b}};s))$ and $\hat{\eta}(s):=f(\emph{\textbf{b}};f(\emph{\textbf{a}};s))$. Then, $\eta(s)$ is the inverse of $\hat{\eta}(s)$. Let $\eta_m(s)$ and $\hat{\eta}_m(s)$ be the $m$'th iterate of $\eta(s)$ and $\hat{\eta}(s)$. Then, $\eta_m(\hat{\eta}_m(s))=s,~\eta_{m+1}(s)\geq \eta_m(s)$ and for $0<s<1$, $\eta_m(s)\uparrow 1$ and $\hat{\eta}_m(s)\downarrow 0$. Fix $0<t_0<1$. Then $t_m:=\eta_m(t_0)\uparrow 1$. Also since $\tilde{Q}(s)$ satisfies (\ref{eq3-5}), (\ref{eq3-6}),
\begin{equation*}
\begin{split}
I_m&
=\int_{t_m}^{t_{m+1}}\tilde{Q}(s)k(s)ds
=\int_{t_m}^{t_{m+1}}\frac{Q(f(\emph{\textbf{a}};s))}{a_1}k(s)ds\\&
=\int_{f(\emph{\textbf{a}};t_m)}^{f(\emph{\textbf{a}};t_{m+1})}
Q(u)\frac{k(g(\emph{\textbf{a}};u))g'(\emph{\textbf{a}};u)du}
{a_1}=\int_{t_{m-1}}^{t_m}\tilde{Q}(s)\frac{k(\eta(s))\eta'(s)}
{a_1b_1}ds\\&
=\int_{t_{m-1}}^{t_m}\tilde{Q}(s)k(s)\frac{k(\eta(s))\eta'(s)}
{a_1b_1k(s)}ds.
\end{split}
\end{equation*}
Since $\eta'(s)=1/\hat{\eta}'(s)$ and $|\log s|/(1-s)\to 1$ as $s\uparrow1$,
$$(k(\eta(s))\eta'(s)/(a_1b_1k(s))\to1/(a_1b_1m_a^r m_b^r),
$$
where $m_a=f'(\emph{\textbf{a}};1)$, $m_b=f'(\emph{\textbf{b}};1)$. Thus if $a_1b_1m_a^rm_b^r>1$, then for any $0<(a_1b_1m_a^rm_b^r)^{-1}<\lambda<1$, there exists an $m_1$ such that $k(\eta(s))\eta'(s)/(a_1b_1k(s))<\lambda$ for all $s\geq \eta_{m_1}(t_0)$. Thus, $I_m\leq\lambda I_{m-1}$ for $m\geq m_1+2$. Hence,
\begin{equation*}
\begin{split}
\sum\limits_{m=m_1+2}^\infty I_m\leq I_{m_1+1}\sum\limits_{j=1}^\infty\lambda^j<\infty.
\end{split}
\end{equation*}
Therefore,
$$
\int_0^1\tilde{Q}(s)k(s)ds\leq \int_0^{t_{m_1}}\tilde{Q}(s)k(s)ds+\int_{t_{m_1}}^1\tilde{Q}(s)k(s)
ds<\infty.
$$
The proof is complete. \hfill $\Box$
\end{proof}
\par
\begin{corollary}\label{cor4.1}
 Assume $a_1,b_1>0$ and $\sum\limits_{j=1}^{\infty}j^{2r+\delta}a_j+\sum\limits_{j=1}
 ^{\infty}j^{2r+\delta}b_j<\infty$ for some $r\geq 1$ and $\delta>0$ such that $a_1m_a^r$, $b_1m_b^r>1$. Then $(\ref{eq4-12}), (\ref{eq4-13})$ hold.
\end{corollary}
\par
\begin{proof}\ Since $\sum\limits_{j=1}^{\infty}j^{2r+\delta}a_j
+\sum\limits_{j=1}
 ^{\infty}j^{2r+\delta}b_j<\infty$ for some $r\geq 1$ and $\delta>0$, we know that
\begin{equation*}
C_{r,1}:= \sup\limits_kE|\sqrt{k}\frac{(\overline{X}_k-m_a)}{\sigma_a}|^{2r}
<\infty
\end{equation*}
and
\begin{equation*}
C_{r,2}:= \sup\limits_kE|\sqrt{k}\frac{(\overline{Y}_k-m_b)}{\sigma_b}|^{2r}
<\infty,
\end{equation*}
where $\{\overline{X}_k:k\geq 1\}$ and $\{\overline{Y}_k:k\geq 1\}$ are given in Theorem~\ref{th4.1}. Then by Markov's inequality, we get
\begin{eqnarray*}
\phi_a(k,\varepsilon)\leq\frac{E|\sqrt{k}
(\overline{X}_k-m_a)|^{2r}}{(\varepsilon\sqrt{k})^{2r}}
\leq\frac{C_{r,1}}{\varepsilon^{2r}k^r}.
\end{eqnarray*}
and
\begin{eqnarray*}
\phi_b(k,\varepsilon)\leq\frac{E|\sqrt{k}
(\overline{Y}_k-m_b)|^{2r}}{(\varepsilon\sqrt{k})^{2r}}
\leq\frac{C_{r,2}}{\varepsilon^{2r}k^r}.
\end{eqnarray*}
Let $C_r=max\{C_{r,1},C_{r,2}\}$. Hence, we have $\phi_b(k,\varepsilon),~\phi_a(k,\varepsilon)\leq\frac{C_r}{
\varepsilon^{2r}k^r}$. Applying Theorem~\ref{th4.2} yields (\ref{eq4-12}) and (\ref{eq4-13}). The proof is complete.
\hfill $\Box$
\end{proof}
\par
Now we consider (\ref{eq1-5}), i.e., the long-time behaviour of $W_n$. Let $\{\tilde{Z}_n:n\geq 0\}$ be the $\emph{\textbf{b}}\!*\!\emph{\textbf{a}}$-Galton-Watson system. Define
\begin{eqnarray*}
\tilde{W}_n=\tilde{\Gamma}_n^{-1}\tilde{Z}_n,\ \ n\geq 0,
\end{eqnarray*}
where $\tilde{\Gamma}_n$ is given in (\ref{eq4-7}).
Similar as the argument about $W_n$, $\tilde{W}_n$ is also an integrable martingale and hence converges to some random variable $\tilde{W}$.
\par
\begin{theorem}\label{th4.3}
 Assume that $f(\textbf{a};e^{\theta_0}), f(\textbf{b};e^{\theta_0})<\infty$ for some $\theta_0>0$. Then there exists $\theta_1>0$ such that
\begin{equation}\label{eq4-14}
C_1=\sup\limits_nE[\exp(\theta_1W_n)]<\infty
\end{equation}
and
\begin{equation}\label{eq4-15}
C_2=\sup\limits_nE[\exp(\theta_1\tilde{W}_n)]<\infty.
\end{equation}
\end{theorem}
\par
\begin{proof}\ Since $K:=f(\emph{\textbf{a}};s_0)<\infty$ for $s_0=e^{\theta_0}$, we know that $f_2(\emph{\textbf{ab}};s)\leq K$ if $0\leq f(\emph{\textbf{b}};s)\leq s_0$, that is, if $0\leq s\leq g(\emph{\textbf{b}};s_0)$. Similarly, $f_3(\emph{\textbf{ab}};s)\leq K$ if $0\leq f(\emph{\textbf{a}};s)\leq g(\emph{\textbf{b}};s_0)$, that is, if $0\leq s\leq g(\emph{\textbf{a}};g(\emph{\textbf{b}};s_0))$. More generally,
\begin{equation*}
f_n(\emph{\textbf{ab}};s)\leq K\ \ ~if~\ 0\leq s\leq g_{n-1}(\emph{\textbf{ba}};s_0).
\end{equation*}
Now, since $W_n=Z_n/\Gamma_n$, $E[e^{\theta W_n}|Z_0=1]=f_n(\emph{\textbf{ab}};e^{\theta/\Gamma_n})$.
Thus
$$
E[\exp(\theta W_n)|Z_0=1]\leq K,
$$
if $\theta\leq \Gamma_n\log g_{n-1}(\emph{\textbf{ba}};s_0).$ Since $g_n(\emph{\textbf{ba}};s_0)\downarrow 1$, $\log g_n(\emph{\textbf{ba}};s_0)\sim (g_n(\emph{\textbf{ba}};s_0)-1)$. By Proposition~\ref{prop4.2}, $f(\emph{\textbf{a}};s_0)<\infty$ for $s_0>1$ implies $\Gamma_n\log g_{n-1}(\emph{\textbf{ba}};s_0)\to m_a\tilde{R}(s_0)$, which is positive and finite. Because of $g_n(\emph{\textbf{ba}};s_0)>1$ for all $n\geq1$, we can choose
\begin{equation*}
\theta_1=\inf\limits_n \Gamma_n\log g_{n-1}(\emph{\textbf{ba}};s_0)~~and~~C_1=K.
\end{equation*}
(\ref{eq4-14}) is proved. A similar argument yields (\ref{eq4-15}).
The proof is complete. \hfill $\Box$
\end{proof}
\par
The next result shows that the convergence rate of $P(|W_n-W|>\varepsilon\mid Z_0=1)$ is supergeometric.
\par
\begin{theorem}\label{th4.4} Let $f(\textbf{a};e^{\theta_0})+f(\textbf{b};e^{\theta_0})<\infty$ for some $\theta_0>0$. Then there exist constants $C_4$ and $\lambda>0$ such that
\begin{equation}\label{eq4-16}
P(|W_n-W|>\varepsilon)\leq C_4\exp(-\lambda\varepsilon^{2/3}\Gamma_n^{1/3}).
\end{equation}
\end{theorem}
\par
\begin{proof} First we need two estimates. Denote
\begin{eqnarray*}
\phi(\theta)=E[\exp(\theta W)]\ \ \ and \ \ \ \tilde{\phi}(\theta)=E[\exp(\theta \tilde{W})],
\end{eqnarray*}
which are finite for all $\theta\leq \theta_1$. So, if $\{W^{(i)}\}_1^\infty ,\{\tilde{W}^{(i)}\}_1^\infty$ are $i.i.d.$ copies of $W$ and $\tilde{W}$ respectively, $S_k=\sum_{i=1}^k(W^{(i)}-1)$, $\tilde{S}_k=\sum_{i=1}^k(\tilde{W}^{(i)}-1)$, then for $\theta\leq \theta_1$,
\begin{equation*}
E[\exp(\theta(S_k/\sqrt{k}))]=(\phi(\frac{\theta}{\sqrt{k}})
e^{-\theta/\sqrt{k}})^k=(1+\frac{1}{k}
\frac{(\phi(\theta/\sqrt{k})e^{-\theta/\sqrt{k}}-1)}
{(\theta^2/k)}\theta^2)^k,
\end{equation*}
\begin{equation*}
E[\exp(\theta(\tilde{S}_k/\sqrt{k}))]=(\tilde{\phi}
(\frac{\theta}{\sqrt{k}})e^{-\theta/\sqrt{k}})^k
=(1+\frac{1}{k}\frac{(\tilde{\phi}(\theta/\sqrt{k})
e^{-\theta/\sqrt{k}}-1)}{(\theta^2/k)}\theta^2)^k.
\end{equation*}
However, since
$$
\lim\limits_{u\to 0}(\phi(u)e^{-u}-1)/u^2
=\frac{1}{2}Var(W)<\infty
$$
and
$$
\lim\limits_{u\to0}(\tilde{\phi}(u)e^{-u}-1)/u^2
=\frac{1}{2}Var(\tilde{W})<\infty,
$$
we have
$$
\sup\limits_{|u|\leq1}|(\phi(u)e^{-u}-1)/u^2|=:c_1<\infty
$$
and
$$
\sup\limits_{|u|\leq1}|(\tilde{\phi}(u)e^{-u}-1)/u^2|=:c_2<\infty.
$$
If $\theta_2=\min(\theta_1,1)$, $c=\max(c_1 ,c_2)$, then
$$
\sup\limits_{|\theta|\leq\theta_2}|\phi(\theta/\sqrt{k})
e^{-\theta/\sqrt{k}}|^k,\ \sup\limits_{|\theta|\leq\theta_2}|\tilde{\phi}(\theta/\sqrt{k})
e^{-\theta/\sqrt{k}}|^k\leq e^{c}=:C_3.
$$
Here we have used the fact that for $x>0$, $(1+x/k)^k\leq e^x$.
\par
Now we proceed with the proof the Theorem~\ref{th4.4}. We begin by noting that (see Theorem 2 on page 55 of Arthreya~\cite{AKN72})
\begin{equation*}
\begin{split}
W-W_n&=
\lim\limits_{m\to\infty}(W_{n+m}-W_n)\\&=
\begin{cases}
\Gamma_n^{-1}\sum\limits_{j=1}^{Z_n}(\tilde{W}^{(j)}-1),\ & if\ n\ is\ odd,\\
\Gamma_n^{-1}\sum\limits_{j=1}^{Z_n}(W^{(j)}-1),\ & if\ n\ is\ even,\\
\end{cases}
\end{split}
\end{equation*}
where $W^{(j)}$ (or $\tilde{W}^{(j)}$ if $n$ is odd) is the limit $r.v$ in the line of descent initiated by the $j$'th parent of the $n$'th generation of $\{Z_n\}$. By conditional independence,
\begin{equation*}
P(W-W_n>\varepsilon|Z_0,Z_1,\cdots,Z_n)=
\begin{cases}
\tilde{\psi}(Z_n,\Gamma_n\varepsilon),\ & if\ n\ is\ odd,\\
\psi(Z_n,\Gamma_n\varepsilon),\ & if\ n\ is\ even,\\
\end{cases}
\end{equation*}
where $\psi(k,\eta)=P(S_k\geq\eta) ,\tilde{\psi}(k,\eta)=P(\tilde{S}_k\geq\eta)$. However, by the above estimation,
\begin{equation*}
\begin{split}
P(S_k\geq\eta)
=P(\frac{S_k}{\sqrt{k}}\geq\frac{\eta}{\sqrt{k}})
\leq C_3\exp(-\frac{\theta_2\eta}{\sqrt{k}}),
\end{split}
\end{equation*}
and
\begin{equation*}
P(\tilde{S}_k\geq\eta)\leq C_3\exp(-\frac{\theta_2\eta}{\sqrt{k}}).
\end{equation*}
Thus,
\begin{equation*}
\begin{split}
P(W-W_n>\varepsilon)&=
\begin{cases}
E[\tilde{\psi}(Z_n,\Gamma_n\varepsilon)],\ if\ n\ is\ odd\\
E[\psi(Z_n,\Gamma_n\varepsilon)],\ if\ n\ is\ even\\
\end{cases}\\
&\leq
C_3E[\exp(-\frac{\theta_2 \Gamma_n\varepsilon}{\sqrt{Z_n}})]\\
&=
C_3E[\exp(-\theta_2\varepsilon \Gamma_n^{1/2}\frac{1}{\sqrt{W_n}})].
\end{split}
\end{equation*}
For $\lambda>0$,
\begin{equation*}
\begin{split}
E[\exp(-\lambda(1/\sqrt{W_n}))]&
=\lambda\int_0^\infty e^{-\lambda u}P(\frac{1}{\sqrt{W_n}}\leq u)du\\
&=\lambda\int_0^\infty e^{-\lambda u}P(W_n\geq\frac{1}{u^2})du\\&
\leq\lambda C_1\int_0^\infty e^{-\lambda u}\exp(-\frac{\theta_1}{u^2})du\ \ ~(by~Theorem~\ref{th4.3})\\
&=C_1\int_0^\infty e^{-t}\exp(-\frac{\theta_1\lambda^2}{t^2})dt.
\end{split}
\end{equation*}
Thus,
\begin{equation*}
P(W-W_n>\varepsilon)\leq C_3C_1\int_0^\infty e^{-t}\exp(-\frac{\theta_1\lambda_n^2}{t^2})dt,
\end{equation*}
where $\lambda_n=\theta_2\varepsilon \Gamma_n^{1/2}$. However, for $\lambda>0$,
\begin{equation*}
I(\lambda):=\int_0^\infty e^{-t}e^{-\lambda^2/t^2}dt=\int_0^{k(\lambda)}+\int_{k(\lambda)}
^\infty\leq\exp(-\frac{\lambda^2}{k^2(\lambda)})+e^{-k(\lambda)}.
\end{equation*}
Choose $k(\lambda)=\lambda^{2/3}$. Then $I(\lambda)\leq2\exp(-\lambda^{2/3})$. Thus
\begin{equation*}
P(W-W_n>\varepsilon)\leq 2C_3C_1\exp(-(\sqrt{\theta_1}\theta_2\varepsilon \Gamma_n^{1/2})^{2/3})=C_4\exp(-\lambda \Gamma_n^{1/3}\varepsilon^{2/3}),
\end{equation*}
where $C_4=2C_3C_1$, $\lambda=(\sqrt{\theta_1}\theta_2)^{2/3}$. Similar arguments hold for $P(W_n-W>\varepsilon)$. Hence, (\ref{eq4-16}) is proved.
The proof is complete. \hfill $\Box$
\end{proof}
\par
Finally, we consider (\ref{eq1-6}). The next result shows that, conditioned on $W>0$, the convergence rate of $P(|\frac{Z_{n+1}}{Z_n}-D_n|>\varepsilon\mid W\geq \delta)$ is supergeometric.
\par
\begin{theorem}\label{th4.5}
Let $f(\textbf{a};e^{\theta_0})+f(\textbf{b};e^{\theta_0})<\infty$ for some $\theta_0>0$. Then there exist constants $C_5$ and $\lambda>0$ such that for all $\varepsilon>0$, $\delta>0$, we can find $0<I(\varepsilon)<\infty$ such that
\begin{equation}\label{eq4-17}
P(|\frac{Z_{n+1}}{Z_n}-\omega_n|>\varepsilon\mid W\geq \delta)
\leq C_5\exp(-\delta\gamma I(\varepsilon)\Gamma_n)
+C_4\exp(-\lambda(\delta(1-\gamma))^{2/3}\Gamma_n^{1/3})
\end{equation}
for every $0<\gamma<1$ and hence $($\ for $\gamma=1/2$\ $)$ there exists constant $C_6>0$ such that
\begin{eqnarray}\label{eq4-18}
P(|\frac{Z_{n+1}}{Z_n}-\omega_n|>\varepsilon\mid W\geq \delta)\leq C_6\exp(-\lambda(\delta/2)^{2/3}\Gamma_n^{1/3}).
\end{eqnarray}
\end{theorem}
\par
\begin{proof}
It is easy to see that
\begin{eqnarray*}
&& P(|\frac{Z_{n+1}}{Z_n}-\omega_n|>\varepsilon\mid W\geq \delta)\\
&=&P(|\frac{Z_{n+1}}{Z_n}-\omega_n|>\varepsilon,W\geq \delta)\frac{1}{P(W\geq \delta)}\\
&=&p_{\delta}[P(|\frac{Z_{n+1}}{Z_n}-\omega_n|>\varepsilon,W_n\leq \delta\gamma,W\geq \delta)+P(|\frac{Z_{n+1}}{Z_n}-\omega_n|>\varepsilon ,W_n\geq \delta\gamma,W\geq \delta)]\\
&=:&p_{\delta}(\delta_{n1}+\delta_{n2}),
\end{eqnarray*}
where $0<\gamma<1$ and $p_{\delta}=1/P(W\geq \delta)$. Clearly,
\begin{equation*}
\delta_{n2}\leq P(|\frac{Z_{n+1}}{Z_n}-\omega_n|>\varepsilon ,W_n\geq \delta\gamma)
\leq C_5\exp(-\delta\gamma I(\varepsilon)\Gamma_n),
\end{equation*}
where $C_5$ and $I(\varepsilon)$ are such that $P(|\overline{Y}_k|\geq \varepsilon)\leq C_5e^{-kI(\varepsilon)}$, $P(|\overline{X}_k|\geq \varepsilon)\leq C_5e^{-kI(\varepsilon)}$ and $\overline{Y}_k=\frac{\sum_{i=1}^kY_i}{k}$, $\{Y_i\}$ being $i.i.d.$ as $Z_{1}^{\emph{\textbf{b}}}-m_b$, $\overline{X}_k=\frac{\sum_{i=1}^kX_i}{k}$, $\{X_i\}$ being $i.i.d.$ as $Z_{1}^{\emph{\textbf{a}}}-m_a$. (Such $C_5$ and $I(\varepsilon)$ exist by Chernoff type bounds since $f(\emph{\textbf{a}};e^{\theta_1})$, $f(\emph{\textbf{b}};e^{\theta_1})<\infty$ for some $\theta_1>0$). Now
\begin{equation*}
\begin{split}
\delta_{n1}&\leq P(W-W_n\geq \delta(1-\gamma))\\&
\leq C_4\exp(-\lambda(\delta(1-\gamma))^{2/3}\Gamma_n^{1/3})\ ~( by\ Theorem~\ref{th4.4}).
\end{split}
\end{equation*}
Therefore,
\begin{equation*}
P(|\frac{Z_{n+1}}{Z_n}-D_n|>\varepsilon|W\geq \delta)\leq p_{\delta}(C_5\exp(-\delta\gamma I(\varepsilon) \Gamma_n)+C_4\exp(-\lambda(\delta(1-\gamma))^{2/3}\Gamma_n^{1/3}).
\end{equation*}
(\ref{eq4-17}) is proved. Since the only condition on $\gamma$ is that $0<\gamma<1$ and the second term goes to zero slower than the first term, we can say that there exist $C_6$ and $\lambda~(C_6$ may depend on $\gamma$ ) such that
\begin{equation*}
P(|\frac{Z_{n+1}}{Z_n}-D_n|>\varepsilon|W\geq \delta)\leq C_6\exp(-\lambda(\delta(1-\gamma))^{2/3}\Gamma_n^{1/3}).
\end{equation*}
Taking $\gamma=\frac{1}{2}$ yields (\ref{eq4-18}). The proof is complete. \hfill $\Box$
\end{proof}
\par
\noindent
{\bf Acknowledgement}\ \ This work is supported by the National Natural Science Foundation of China (No. 11771452, No. 11971486).
\par
\noindent
{\bf Declarations}
\hspace*{0.5cm}
\par
\noindent
{\bf  Ethics Approval}\ Not applicable.
\par
\noindent
{\bf Conflict of Interests}\ The authors declare that they have no conflict of interest.
\par
\noindent
{\bf Data Availability}\ Not applicable.
\par
\noindent
{\bf Funding}\ Funding provided by the National Natural Science Foundation of China.

\par

\end{document}